\documentclass[preprint,12pt]{elsarticle}

\usepackage{fullpage}

\usepackage[latin1]{inputenc}
\usepackage{amsmath,amssymb,amsthm}
\usepackage{graphicx,graphics,color,epsfig}
\usepackage{hyperref}

\newcommand{\CN}{\mathcal N}

\newcommand{\CS}{\mathcal S}

\newcommand{\CU}{\mathcal U}
\newcommand{\CV}{\mathcal V}
\newcommand{\CW}{\mathcal W}

\newcommand{\BE}{\mathbb E}

\newcommand{\BN}{\mathbb N}

\newcommand{\BR}{\mathbb R}

\theoremstyle{plain}
\newtheorem{theorem}{Theorem}

\theoremstyle{definition}

\theoremstyle{remark}

\newtheorem{example}{Example}

\graphicspath{{../figures/}}

\begin{document}

\begin{frontmatter}

\title{Piecewise linear car-following modeling}


\author{Nadir Farhi}
\ead{nadir.farhi@ifsttar.fr}


\date{}

\address{Universit\'e Paris-Est, IFSTTAR, GRETTIA, F-93166 Noisy-le-Grand, France.}

\begin{abstract}
We present a traffic model that extends the linear car-following model as well as the min-plus traffic
model (a model based on the min-plus algebra).
A discrete-time car-dynamics describing the traffic on a 1-lane road without passing is
interpreted as a dynamic programming equation of a stochastic optimal control problem of a Markov chain.
This variational formulation permits to characterize the stability of the car-dynamics and to calculte
the stationary regimes when they exist. The model is based on a piecewise linear approximation of the fundamental
traffic diagram.
\end{abstract}


\end{frontmatter}

\section{Introduction}

Car-following models are microscopic traffic models that describe the car-dynamics with stimulus-response
equations expressing the drivers' behavior.
Each driver responds, by choosing its speed or acceleration, to a given stimulus
that can be composed of many factors such as inter-vehicular distances, relative velocities, instantaneous
velocities, etc.
We present in this article a car-following model that extends the linear car-following
model~\cite{CHM58,HMPR59,GHP59,GHR61}, as well as the min-plus traffic model~\cite{LMQ05}.

The vehicular traffic on a 1-lane road without passing is described by
discrete-time dynamics, which are interpreted as dynamic programming equations associated to stochastic optimal
control problems of Markov chains. The discrete-time variational formulation we make here is similar to the
time-continuous one used by Daganzo and Geroliminis~\cite{DG08} to show the existence of a concave macroscopic
fundamental diagram.

We are concerned in this article by microscopic traffic modeling with a Lagrangian description of the traffic
dynamics, where the function $x(n,t)$, giving the position of car $n$ at time $t$ (or the
cumulated distance traveled by a car $n$ up to time $t$), is used.
In the macroscopic kinematic traffic modeling, Eulerian descriptions of the traffic dynamics are usually used,
with the function $n(t,x)$ giving the cumulated number of cars passing through position $x$ up to
time $t$ (which coincides with the Moscowitz function~\cite{Dag06} in the case of traffic without passing).
The combination of a conservation law with an equilibrium law gives the well known first order traffic model
of Lighthill, Whitham and Richards~\cite{LW55,Ric56}.

In this introduction, we first notice that the same approach used in macroscopic traffic
modeling, combining a conservation law with an equilibrium law, can be used to derive microscopic traffic models.
By this, we introduce car-following models and give a review on basic and most important existing ones.
In particular, we recall the linear and the min-plus car-following models, which are
particular cases of the model we present here. Finally, we recall a theoretical result on non-expansive
and connected dynamic systems which we need in our developments, and give the outline of this article.

The first order partial derivative of the function $x(n,t)$ in time, denoted $\dot{x}(n,t)$
expresses the velocity $v(n,t)$ of car $n$ at time $t$. The first order differentiation of
$x(n,t)$ in the car numbers $\Delta x(n,t)$ expresses the inverse of the
inter-vehicular distance $y(n,t) = x(n-1,t)-x(n,t)$ between cars $n$ and $n-1$ at time $t$.
The equality of the second derivatives $\Delta \dot{x}(n,t)$ and
$\dot{\Delta} x(n,t)$ gives then the following conservation law of distance.
\begin{equation}\label{conservlaw2}
  \dot{y}(n,t)= - \Delta v(n,t).
\end{equation}

If we assume that a fundamental diagram $V_e$, giving the velocity $v$ as a function
of the inter-vehicular distance $y$ ($v=V_e(y)$) at the stationary traffic, exists,
and that the diagram $V_e$ holds also on the transient traffic, then we have
\begin{equation}\label{doty}
  \dot{y}(n,t)=\dot{v}(n,t) / V_e'(y),
\end{equation}
where $V_e'(y)$ denotes the derivative of $V_e$ with respect to $y$.

By combining~(\ref{conservlaw2}) and~(\ref{doty}) we obtain the model:
\begin{equation}\label{carf1}
   \dot{v}(n,t) = - V_e'(y(n,t)) \Delta v(n,t).
\end{equation}

(\ref{carf1}) is a car-following model that gives the acceleration of car $n$ at time $t$ as a response to a
stimulus composed of the relative speed $\Delta v(n,t)$ and the term $ V_e'(y(n,t))$.
For example, if $V_e(y)=v_0 \exp{-a/y}$, where $v_0$ denotes the free (or desired) velocity,
and $a$ is a parameter, then $V'(y)=a V(y)/y^2$ and~(\ref{carf1}) gives a particular case of the Gazis, Herman, and Rothery model~\citep{GHR61}.

The simplest car-following model is the linear one, where the car dynamics are written
\begin{equation}\label{Lin2}
  \dot{x}_n(t+T)=a (x_{n-1}(t)-x_n(t))+b,
\end{equation}
where $T$ is the reaction time, $a$ is a sensitivity parameter, and $b$ is a constant.
That model can be derived from~(\ref{carf1}) with a linear fundamental diagram $V_e$.
The stability of the linear car-following model~(\ref{Lin2}) and the existence of a stationary
regime have been treated in~\cite{HMPR59}.

Almost all car following models are based on the assumption of the existence of a behavioral law $V_e$.
The latter has been taken linear in~\cite{Green35,HMPR59}, logarithmic in~\cite{Gre59}, exponential in~\cite{New61},
and with more complicated forms in other works.
Bando et al.~\cite{BHNSY95} used the sigmiudal function
\begin{equation}\label{Banlaw}
  V_e(y)=tanh(y-h)+tanh(h),
\end{equation}
where $tanh$ denotes the hyperbolic tangent function, and $h$ is a constant.

Kerner and Konh\"auser\citep{KK93}, Hermann and Kerner~\citep{HK98}, and then Lenz et al.\citep{LWS99}, and Hoogendoorn et al.~\citep{HOS06}
have used the function
\begin{equation}\label{Kerlaw}
  V_e(y)=v_0 \left\{\left[1+exp\left(\frac{1,000}{\gamma \cdot y}-\frac{10}{2.1}\right)\right]^{-1}-5.34\cdot 10^{-9}\right\},
\end{equation}
where free velocity $v_0$ and $\gamma$ are parameters estimated from data. In~\citep{LWS99}, $\gamma$ is taken equal to $7.5$.
See also~\cite{Hel59,Hel61,THH00}.

The min-plus traffic model~\cite{LMQ05} is a microscopic discrete-time car-following model
based on the min-plus algebra~\cite{BCOQ92}. It consists in the following dynamics.
\begin{equation}\label{minplus}
  x_n(t+1)=\min\{x_n(t)+v_0, x_{n-1}(t)-\sigma\},
\end{equation}
where $v_0$ is the free velocity and $\sigma$ is a safety distance.

The idea of the model~(\ref{minplus}) is that the car dynamics is linear in the min-plus algebra~\cite{BCOQ92},
where the addition is the operation ``$\min$'' and the product is the standard addition ``$+$''.
In~\citep{LMQ05}, the dynamic system~(\ref{minplus}) is written in
min-plus notations $x(t)=A\otimes x(t-1)$, where $A$ is a min-plus matrix
and $\otimes$ is the min-plus product of matrices.
It is then proved that the average growth rate vector per
time unit of the system, defined by $\chi=\lim_{t\to\infty}x(t)/t$ satisfies
$\chi=(\bar{v},\bar{v},\cdots,\bar{v})$, where $\bar{v}$
is the unique min-plus eigenvalue of the matrix $A$, interpreted as the stationary car-speed.
The fundamental diagram is then obtained
\begin{align}
  & \bar{v}=\min(v_0, \bar{y}-\sigma) \label{mplaw1} \\
  & \bar{q}=\min(v_0 \bar{\rho}, 1-\sigma \bar{\rho}), \label{mplaw2}
\end{align}
where $\bar{y}, \bar{q}$ and $\bar{\rho}$ denote respectively the equilibrium inter-vehicular distance,
car-flow and car-density.

In~(\ref{minplus}), for $t\geq 0$, if $x_{n-1}(t)-x_n(t) > v_0+\sigma$, then the dynamics
is $x_n(t+1)-x_n(t)=v_0$. If $x_{n-1}(t)-x_n(t) \leq v_0+\sigma$ then the dynamics is
$x_n(t+1)-x_n(t)=x_{n-1}(t)-x_n(t) -\sigma$. Therefore,~(\ref{minplus}) is linear in both
phases of free and congested traffic.
The min-plus model~(\ref{minplus}) permits to distinguish two phases in which the traffic
dynamics is linear, but with a sensitivity parameter ($a$ in~(\ref{Lin2})) equals to $1$
for each phase. The linear model~(\ref{Lin2}) is not constrained in the sensitivity parameter
value, but it permits the modeling of only one traffic phase.
The model we present here extends~(\ref{Lin2}) and~(\ref{minplus}) in a way that
an arbitrary number of traffic phases can be modeled, with flexibility in the sensitivity
parameter value on each phase.

For our model, we apply a similar but more general approach than the min-plus one used
to analyze the dynamic system~(\ref{minplus}).
Indeed, the dynamics~(\ref{minplus}) is additive homogeneous of
degree one~\footnote{A dynamic system $x(t)=f(x(t-1))$ is additive homogeneous
of degree 1 if $f$ is so, that is if $\forall x\in\BR^n, \forall \lambda \in\BR,
f(\lambda+x)=\lambda+f(x)$.} and is monotone~\footnote{A dynamic system
$x(t)=f(x(t-1))$ is monotone if $f$ is so, that is if
$\forall x_1, x_2 \in\BR^n, x_1\leq x_2 \Rightarrow f(x_1)\leq f(x_2)$, where the order $\leq$
is pointwise in $\BR^n$.}. It is then non expansive~\footnote{A dynamic system $x(t)=f(x(t-1))$ is
non expansive if $f$ is so, that is if there exists a norm $||\cdot||$ in $\BR^n$ such that
$\forall x_1,x_2\in\BR^n, ||f(x_2)-f(x_1)||\leq ||x_2-x_1||$.}~\citep{CT80}.
The stability of the dynamic system~(\ref{minplus}) is thus guaranteed from its
non expansiveness. Moreover,~(\ref{minplus}) is
connected (or communicating)~\footnote{An additive homogeneous of degree 1 and monotone dynamic
system $x(t)=f(x(t-1))$ with $x\in\BR^n$ is connected if its associated graph is strongly connected.
The graph associated to that dynamic system is the graph with $n$ nodes and whose arcs
are determined as follows. There exists an arc from a node $i$ to a node
$j$ if $lim_{\nu\to\infty}f_j(\nu e_i)=\infty$, where $e_i$
denotes the $i^{th}$ vector of the canonic basis of $\BR^n$.}~\citep{GG98a,GG99}.
An important result from~\citep{GK95,GG98a} (Theorem~\ref{theoG} below) permits the analysis of
non expansive and connected dynamic systems.
\begin{theorem}\cite{GK95,GG98a}\label{theoG}
  If a dynamic system $x(t)=f(x(t-1))$ is non expansive and connected, then
  the additive eigenvalue problem $\bar{v}+x=f(x)$ admits a solution $(\bar{v},x)$, where $x$ is
  defined up to an additive constant, not necessarily in a unique
  way, and $\bar{v}\in\BR^n$ is unique. Moreover, the dynamic system admits an average growth
  rate vector $\chi$, which is unique (independent of the initial condition) and given by
  $\chi(f)={}^t(\bar{v},\bar{v},\cdots,\bar{v})$.
\end{theorem}

The model treated in this article can be seen as an extension of the min-plus
model~(\ref{minplus}).
In Section 2, we present the model. It is called \emph{piecewise
linear car-following model} because it is based on a piecewise linear fundamental diagram $V_e$.
The model describes the traffic of cars on a ring road of one lane without passing.
The stability conditions of the dynamic system describing the traffic are determined.
Under those conditions, the car-dynamics are interpreted as a dynamic programming equation~(DPE)
associated to a stochastic optimal control problem of a Markov chain. The DPE is solved
analytically. We show that the individual behavior law $V_e$ supposed
in the model is realized on the collective stationary regime.
Finally, the effect of the stability condition on the shape of the fundamental diagram is shown.

In section 3, we give equivalent results of those given in section 2 for the traffic on an
``open'' road (a highway stretch for example), and conclude with an example, where we simulate
the transient traffic, basing on a piecewise linear approximation of the diagram~(\ref{Kerlaw}).

In~\ref{appB}, we give more details on the duality in traffic modeling, of using the functions
$n(t,x), x(n,t)$ and $t(n,x)$, where $t(n,x)$ denotes the time of passage of the
$n$th car by position $x$.

\section{Piecewise linear car following model}

The behavioral law $V_e$ is an increasing curve bounded by the free speed $v_0$.
Moreover, $V_e(y)=0$ for $y\in[0,y_j]$ where $y_j$ denotes the jam inter-vehicular
distance. We propose here to approximate the curve $V_e$ with a piecewise-linear
curve
\begin{equation}\label{approx}
  V_e(y)=\min_{u\in\CU}\max_{w\in\CW} \{\alpha_{uw} y + \beta_{uw}\},
\end{equation}
where $\CU$ and $\CW$ are two finite sets of indices.
Since $V_e$ is increasing, we have $\alpha_{uw} \geq 0, \forall (u,w)\in\CU\times\CW$.

We are interested here on the discrete-time first-order dynamics
\begin{equation}\label{gen1}
  x_n(t+1)=x_n(t)+\min_{u\in\CU}\{\alpha_u (x_{n-1}(t)-x_n(t))+\beta_u\},
\end{equation}
and
\begin{equation}\label{gen2}
  x_n(t+1)=x_n(t)+\min_{u\in\CU}\max_{w\in\CW}\{\alpha_{uw} (x_{n-1}(t)-x_n(t))+\beta_{uw}\}.
\end{equation}
It is clear that~(\ref{gen2}) extends~(\ref{gen1}). The model~(\ref{gen2}) is also an extension
of both linear model~(\ref{Lin2}) and min-plus model~(\ref{minplus}).

In this article, we characterize the stability of the dynamics~(\ref{gen2}), calculate the stationary regimes,
show that the fundamental diagrams are effectively realized at the stationary regime, and analyze the transient traffic.
We will distinguish two cases: Traffic on a ring road and traffic on an ``open'' road.

\subsection{Traffic on a ring road}

We follow here the modeling of~\citep{LMQ05}.
Let us consider $\nu$ cars moving a one-lane ring road in one direction without passing.
We assume that the cars have the same length that we take here as the unity of distance.
The road is assumed to be of size $\mu$; that is, it can contain at most $\mu$ cars.
The car density on the road is thus $\rho=\nu/\mu$.

\subsection*{Stochastic optimal control model}

We consider here the car dynamics~(\ref{gen1}). That is to say that each car $n$ maximizes its velocity
at time $t$ under the constraints
\begin{equation}\label{cons1}
  x_n(t+1)\leq x_n(t)+ \alpha_u (x_{n-1}(t)-x_n(t))+\beta_u, \quad \forall u\in\CU.
\end{equation}
Each constraint of~(\ref{cons1}) bounds the velocity $x_n(t+1)-x_n(t)$ by a sum of a fixed
term $\beta_u$ and a term depending linearly on the inter-vehicular distance.

Let us first notice that~(\ref{gen1}), on the ring road, is written
\begin{align}
  & x_n(t+1)=x_n(t)+\min_{u\in\CU}\{\alpha_u (x_{n-1}(t)-x_n(t))+\beta_u\}, \quad \text{ for } n\geq 2, \nonumber \\
  & x_1(t+1)=x_1(t)+\min_{u\in\CU}\{\alpha_u (x_{n}(t)+\mu-x_1(t))+\beta_u\}, \nonumber
\end{align}
which can also be written
\begin{align}
  & x_n(t+1)= \min_{u\in\CU}\{(1-\alpha_u) x_n(t)+ \alpha_u x_{n-1}(t)+\beta_u\}, \quad \text{ for } n\geq 2, \nonumber \\
  & x_1(t+1)= \min_{u\in\CU}\{(1-\alpha_u) x_1(t)+ \alpha_u x_{n}(t)+\alpha_u \nu/\rho+\beta_u\}. \nonumber
\end{align}

Let us denote by $M^u, u\in\CU$ the family of matrices defined by
$$M^u=\begin{bmatrix}
        1-\alpha_u & 0 & \cdots & \alpha_u\\
        \alpha_u & 1-\alpha_u & & 0\\
        \vdots & \ddots & \ddots & \\
        0 & 0 & \alpha_u & 1-\alpha_u
      \end{bmatrix},$$
and by $c^u, u\in\CU$, the family of vectors defined by
$$c^u={}^t[\alpha_u\nu/\rho+\beta_u,\;\; \beta_u,\; \cdots,\; \beta_u].$$

The dynamics~(\ref{gen1}) are then written~:
\begin{equation}  \label{gener1}
  x_n(t+1) = \min_{u\in\CU}\{[M^u x(t)]_n+c^u_n\},\quad 1\leq n\leq \nu\;.
\end{equation}

The system~(\ref{gener1}) is additive homogeneous of degree~1 by the definition of the matrices $M^u, u\in\CU$.
It is monotone under the condition that all the components of $M^u, u\in\CU$ are non negative, which is equivalent
to $\alpha_u\in [0,1], \forall u\in\CU$. Hence, under that condition, the system~(\ref{gener1}) is
non expansive.

Moreover, the matrices $M^u, u\in\CU$ are
stochastic~\footnote{We mean here $M^u_{ij}\geq 0, \forall i,j$ and $\sum_j M^u_{ij}=1, \forall i$.}.
Those matrices can then be seen as transition matrices of a controlled Markov chain, where
the set of controls is $\CU$. The connectedness of the system~(\ref{gener1}),
as defined in the previous section, is related to the irreducibility of the Markov chain with transition
matrices $M^u, u\in\CU$. It is easy to check that~(\ref{gener1}) is connected if and only if
$\exists u\in\CU, \alpha_u\in (0,1]$; see~\ref{appA} for the proof.
That condition is interpreted in term of traffic by saying that every car moves by taking into account the
position of the car ahead.

Consequently, under the condition $\forall u\in\CU, \alpha_u\in [0,1]$ and $\exists u\in\CU, \alpha_u\in (0,1]$,
the dynamic system~(\ref{gener1}) is non expansive and is connected.
Therefore, by Theorem~\ref{theoG}, we conclude that the additive eigenvalue problem
\begin{equation}\label{gen-eig}
  \bar{v}+x_n = \min_{u\in\CU}\{[M^u x]_n+c^u_n\},\quad 1\leq n\leq \nu\;,
\end{equation}
describing the stationary regime of the dynamic system~(\ref{gener1}), admits a solution
$(\bar{v},x)$, where $x$ is defined up to an additive constant, not necessarily in a unique
way, and $\bar{v}\in\BR^n$ is unique. Moreover, the dynamic system admits a unique average growth
rate per time unit $\chi$, whose components are all equal and coincide with $\bar{v}$.

The average growth rate per time unit $\chi$ of the system~(\ref{gener1}) is interpreted in term of
traffic as the stationary car-velocity. The additive eigenvector $x$ gives the asymptotic distribution
of cars on the ring. $x$ is given up to an additive constant, since the car-dynamics~(\ref{gener1}) is
additive homogeneous of degree 1. That is to say that if $(\bar{v},x)$ is a solution of~(\ref{gen-eig})
then $(\bar{v}, e+x)$ is also a solution for~(\ref{gen-eig}), for every constant $e\in\BR$.

Let us now give an interpretation of the model in term of ergodic stochastic optimal control.
Indeed, (\ref{gen-eig}) can be seen as a dynamic programming equation of an
ergodic stochastic optimal control problem of a Markov chain with transition matrices
$M^u, u\in\CU$ and costs $c^u, u\in\CU$, and with a set of states $\CN=\{1,2,\cdots,\nu\}$.
The stochastic optimal control problem of the chain is written
\begin{equation}\label{socp1}
  \min_{s\in\CS}\BE \left\{\lim_{T\to +\infty} \frac{1}{T} \sum_{t=0}^{T} c^{u_t}_{n_t}\right\},
\end{equation}
where $\CS$ is a set of feedback strategies on $\CN$. A strategy $s\in\CS$ associates to every
state $n\in\CN$ a control $u\in\CU$ (that is $u_t=s(n_t)$).

The following result gives one solution $(\bar{v},x)$ for the dynamic programming equation~(\ref{gen-eig}).
\begin{theorem}\label{th-gen-2}
  The system~(\ref{gen-eig}) admits a solution $(\bar{v},x)$ given by:
  \begin{align}
    \bar{v} = & \min_{u\in\CU}\{\alpha_u \bar{y}+\beta_u\}, \nonumber \\
    x   = & {}^t[0 \quad \bar{y} \quad 2\bar{y} \quad \cdots \quad (\nu-1)\bar{y}]. \nonumber
  \end{align}
\end{theorem}
\proof First, because of the symmetry of the system~(\ref{gen-eig}), it is natural
  that the asymptotic car-positions $x_n, 1\leq n\leq \nu$ are uniformly distributed on the ring,
  and that the optimal strategy is independent of the state $x$. Let us prove it.

  Let $\bar{u}\in\CU$ be defined by
  $\alpha_{\bar{u}}\bar{y}+\beta_{\bar{u}}=\min_{u\in\CU}\{\alpha_u\bar{y}+\beta_u\}=\bar{v}$.
  Let $x$ be the vector given in Theorem~\ref{th-gen-2}. Then $\forall n\in\{1,2,\cdots,\nu\}$ we have
  $$[M^{\bar{u}} x]_n + c^{\bar{u}}_n = (\alpha_{\bar{u}}\bar{y}+\beta_{\bar{u}}) + x_n
                  = \min_{u\in\CU}(\alpha_{u}\bar{y}+\beta_{u}) + x_n
                  = \min_{u\in\CU} \{[M^{u} x]_n + c^{u}_n\}
                  = \bar{v} + x_n.$$
\endproof

In term of traffic, Theorem~\ref{th-gen-2} shows that the car-dynamics is stable under the condition $\alpha_u\in[0,1]$,
and the average car speed is given by the additive eigenvalue of the asymptotic dynamics in the case where the system is connected.
Moreover, it affirms that the fundamental diagram supposed in the model is realized at the stationary regime.
\begin{align}
  & \bar{v} = \min_{u\in\CU}\{\alpha_u \bar{y}+\beta_u\}, \label{diagv} \\
  & \bar{q} = \min_{u\in\CU}\{\alpha_u+\beta_u \bar{\rho}\}. \label{diagq}
\end{align}

It is important to note here that, up to the assumption $\alpha_u\in [0,1], \forall u\in\CU$,
every concave curve $V_e$ or $Q_e$ can be approximated with~(\ref{diagv}) or~(\ref{diagq}).
Indeed, approximating fundamental diagrams using those formulas is nothing but computing
Fenchel transforms; see ~\cite{Dag06,ABS08}. More precisely, if we denote by
$\CV$ the set $\CV=\{\beta_u, \; u\in\CU\}$ and define the function $g$ by:
$$\begin{array}{llll}
  g: & \CV & \to & \BR\\
     & v=\beta_u & \mapsto & - \alpha_u\;,
\end{array}$$
then
$$q=Q_e(\rho)=\min_{v\in\CV}(\rho v - g(v)) = g^*(\rho),$$
where $g^*$ denotes the Fenchel transform of $g$.

Finally, we note that the min-plus linear model is a particular case of the model presented in
this section, where $\CU=\{u_1,u_2\}$ with $(\alpha_1,\beta_1)=(0,v)$ and
$(\alpha_2,\beta_2)=(1,-\sigma)$.
In this case, the fundamental traffic diagram is approximated with a piecewise linear curve with two segments.

\subsection*{Stochastic game model}

We consider in this section the car dynamics~(\ref{gen2}), again with the assumption
$\forall (u,w)\in\CU\times\CW, \alpha_{uw}\in[0,1]$ and $\exists (u,w)\in\CU\times\CW, \alpha_{uw}\in (0,1]$.
The dynamic system~(\ref{gen2}) is interpreted as a stochastic dynamic programming equation
associated to a stochastic game problem on a controlled Markov chain.
As above, a generalized eigenvalue problem is solved.
The extension we make here approximates non concave fundamental diagrams.

In term of traffic, we take into account the drivers' behavior changing from low densities to high ones.
The difference between these two situations is that in low densities, drivers, moving, or \emph{being able to move} with
high velocities, they try to leave large safety distances between each other, so the safety distances are
maximized; whilst in high densities, drivers, moving, or \emph{having to move} with low velocities, they try to leave
small safety distances between each other in order to avoid jams; so they minimize safety distances.

To illustrate this idea, let us consider the following two dynamics of a given car $n$.
\begin{align}
  & x_n(t+1)= \min\{ x_n(t)+v, x_{n-1}(t)-\sigma\}, \label{dyn1} \\
  & x_n(t+1)= \min\{ x_n(t)+v, \max\{x_{n-1}(t)-\sigma,(x_n(t)+x_{n-1}(t))/2\}\}. \label{dyn2}
\end{align}
The dynamics~(\ref{dyn1}) is a min-plus dynamics which grossly tell that cars move with
their desired velocity $v$ at the fluid regime and they keep a safety distance $\sigma$
at the congested regime. The dynamics~(\ref{dyn2}) distinguishes two situations at the
congested regime:
\begin{itemize}
  \item In a relatively low density situation where the cars are separated by a distance that
    equals at least to $2\sigma$ we have
    $$\max\{x_{n-1}(t)-\sigma,(x_n(t)+x_{n-1}(t))/2\}=x_{n-1}(t)-\sigma.$$
  \item In a high density situation, where the cars are separated by distances less than $2\sigma$
    we have
    $$\max\{x_{n-1}(t)-\sigma,(x_n(t)+x_{n-1}(t))/2\}=(x_n(t)+x_{n-1}(t))/2.$$
    In this case, we accept the cars moving closer but by reducing the approach speed in order to
    avoid collisions. This is realistic.
\end{itemize}

The situation we have considered in~(\ref{dyn2}) is realistic and very simple, but, it cannot
be obtained without introducing a \emph{maximum} operator in the dynamics (i.e. with only minimum operators).
Indeed, with only minimum operators the approach is mechanically reduced with the increasing
of the car-density (in fact this is the concaveness of the fundamental diagram).
Because of the realness of such scenarios, we think that the fundamental traffic diagram should
be composed of two parts, a concave part at the fluid regime, and a convex part at the
congested regime. The dynamics~(\ref{gen2}) generalizes this idea.

The dynamics~(\ref{gen2}) can be written
\begin{equation}\label{gener2}
  x_n(t+1)=\min_{u\in\CU}\max_{w\in\CW}\{[M^{uw} x(t)]_n+c^{uw}_n\},\quad 1\leq n\leq \nu\;,
\end{equation}
where
$$M^{uw}=\begin{bmatrix}
        1-\alpha_{uw} & 0 & \cdots & \alpha_{uw}\\
        \alpha_{uw} & 1-\alpha_{uw} & & 0\\
        \vdots & \ddots & \ddots & \\
        0 & 0 & \alpha_{uw} & 1-\alpha_{uw}
      \end{bmatrix},$$
and
$$c^{u}={}^t[\alpha_{uw}\nu/\rho+\beta_{uw},\;\; \beta_{uw},\; \cdots,\; \beta_{uw}].$$

Similarly, we can easily check that the dynamic system~(\ref{gener2}) is non expansive under the
condition $\alpha_{uw}\in [0,1], \forall (u,w)\in\CU\times\CW$. Its stability is thus guaranteed
under that condition.
If, in addition, $\exists (u,w)\in\CU\times\CW, \alpha_{uw}\in (0,1]$, then the system is
connected.
In this case, we get the same results as in Theorem~\ref{th-gen-2}. That is, the eigenvalue
problem
\begin{equation}\label{gen-eig2}
  \bar{v}+x_n=\min_{u\in\CU}\max_{w\in\CW}\{[M^{uw} x]_n+c^{uw}_n\},\quad 1\leq n\leq \nu
\end{equation}
admits a solution $(\bar{v},x)$ given by:
\begin{align}
  \bar{v} = & \min_{u\in\CU}\max_{w\in\CW}\{\alpha_{uw} \bar{y}+\beta_{uw}\}, \nonumber \\
  x = & {}^t[0 \quad \bar{y} \quad 2\bar{y} \quad \cdots \quad (\nu-1)\bar{y}]. \nonumber
\end{align}
Moreover, the dynamic system~(\ref{gener2}) admits a unique average growth rate vector $\chi$,
whose components are all equal to $\bar{v}$.

The car-dynamics is then stable under the condition $\forall (u,w)\in \CU\times\CW, \alpha_{uw}\in[0,1]$,
and the average car speed is given by the additive eigenvalue of the asymptotic dynamics in the case where the
system is connected (that is, if $\exists (u,w)\in\CU\times\CW, \alpha_{uw}\in (0,1]$).
The behavior law supposed in the model is realized at the stationary regime.
\begin{align}
  & \bar{v} = \min_{u\in\CU}\max_{w\in\CW}\{\alpha_{uw} \bar{y}+\beta_{uw}\}, \label{diagv2} \\
  & \bar{q} = \min_{u\in\CU}\max_{w\in\CW}\{\alpha_{uw}+\beta_{uw} \bar{\rho}\}. \label{diagq2}
\end{align}

In term of stochastic optimal control, the system~(\ref{gen-eig2}) can be seen as a dynamic programming
equation associated to a stochastic game, with two players, on a Markov chain.
The set of states of the chain is again $\CN=\{1,2,\cdots,\nu\}$. The chain is controlled by two players, a minimizer
one with a finite set $\CU$ of controls, and a maximizer one with a finite set $\CW$ of controls.
The transitions and the costs of the chain are given by the matrices $M^{uw}$ and the vectors
$c^{uw}, (u,w)\in\CU\times\CW$ defined above.

The stochastic optimal control problem is
\begin{equation}\label{gamepbm}
  \min\max|_{s\in\CS}\BE\left\{\lim_{T\to +\infty}\frac{1}{T}\sum_{t=0}^{T}c^{u_t w_t}_{n(t)}\right\},
\end{equation}
where $\CS$ is the set of strategies assoicating to every state $n\in\CN$ a couple of commands $(u,w)\in\CU\times\CW$.
It is assumed here that the maximizer knows at each step the decision of the minimizer.

\begin{figure}[htbp]
  \begin{center}
    \includegraphics[width=8cm]{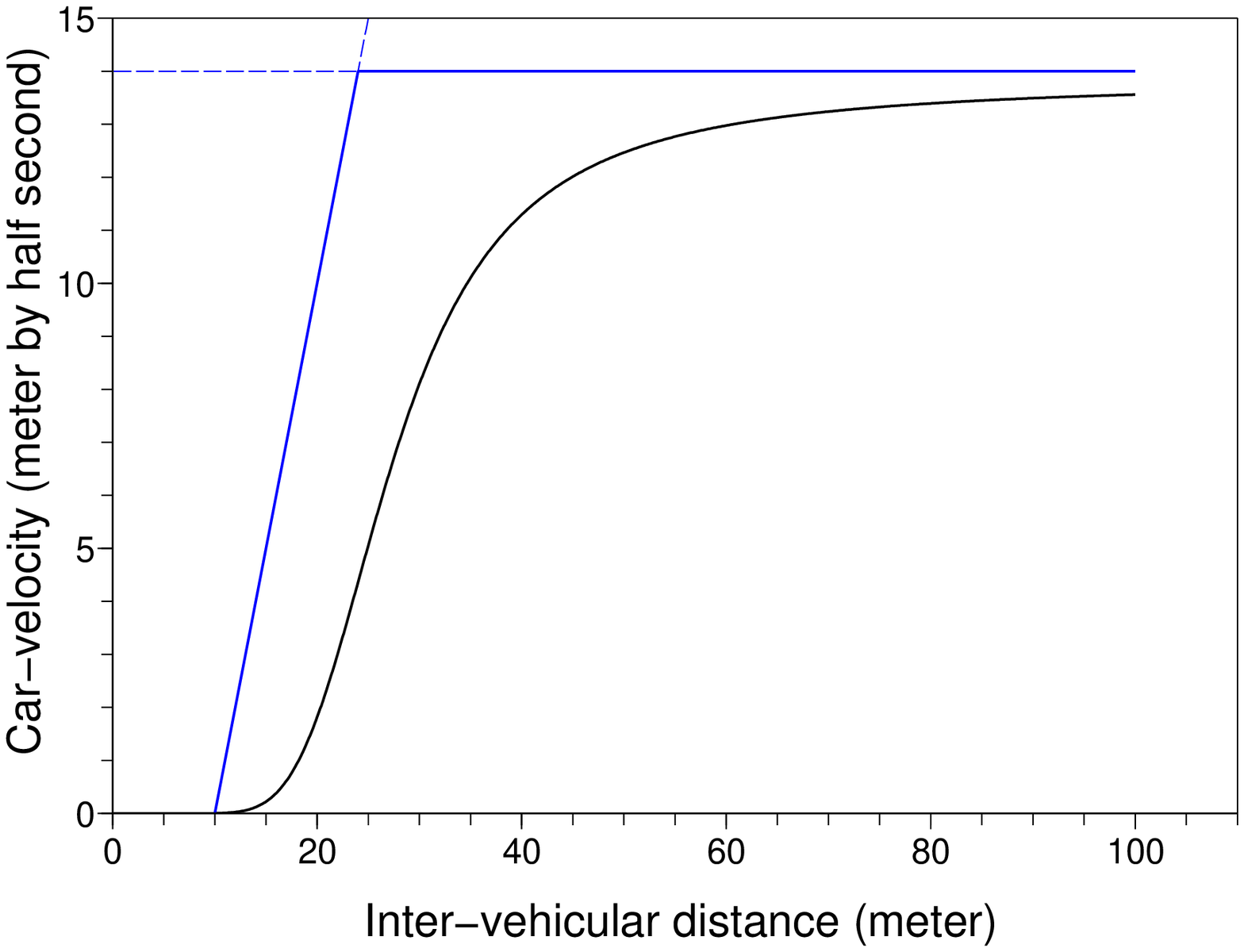}
    \includegraphics[width=8cm]{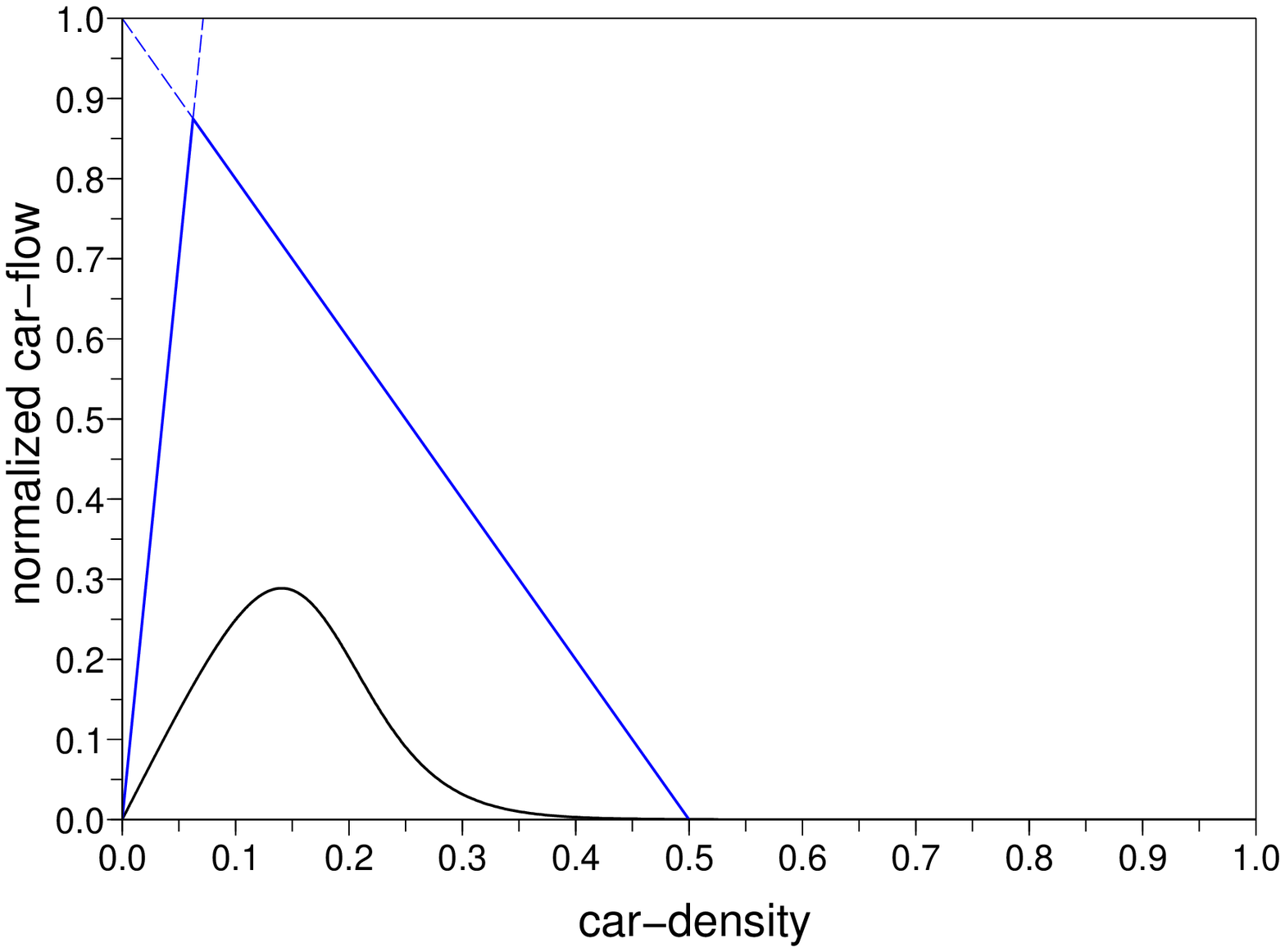}
    \caption{The effect of the stability condition $\alpha_{uw}\in [0,1]$ on the shape of the fundamental diagram.}
    \label{half}
  \end{center}
\end{figure}

We now give a consequence of the stability condition $\alpha_{uw}\in[0,1], \forall (u,w)\in\CU\times\CW$,
on the shape of the fundamental diagrams~(\ref{diagv2}) and~(\ref{diagq2}). As shown in~Figure~\ref{half},
where we have drawn the fundamental diagram~(\ref{Kerlaw}) (with $v_0=14$ meter by half second, and $\gamma=7.5$),
the stability condition puts the curves~(\ref{diagv2}) and~(\ref{diagq2}) in specific respective regions in the plan.
Indeed, for the diagram~(\ref{diagv2}), if we assume that $V_e$ is bounded by $v_0$, $V_e(y)=0, \forall y\in [0,y_j]$, and
that $V_e$ is continuous (and increasing), then starting by the point $(y_j,0)$, one cannot join any point above
the line passing by $(y_j,0)$ and having the slope $1$, with any sequence of segments of slopes $\alpha_{uw}\in[0,1]$.
We can write
$$V_e(y)\leq \max(0, \min(v_0, y - y_0)).$$
Similarly, on the diagram~(\ref{diagq2}), if we assume that $Q_e$ is continuous and $Q_e(\rho)=0, \forall \rho\in [\rho_j,1]$, then
going back from the point $(\rho_j,0)$, one cannot attain any point above the line passing by $(\rho_j,0)$
and $(0,1)$, with a sequence of segments having their ordinates at the origin ($\alpha_{uw})$ in $[0,1]$.
We can write
$$Q_e(\rho) \leq \max(0, \min(v_0 \rho, 1 - \rho/\rho_j)).$$

\section{Traffic on an open road}

We study in this section the traffic on an open road with one lane and without passing.
We are interested in the following dynamics.
\begin{equation}\label{open1}
  \begin{array}{l}
    x_1(t+1)=x_1(t)+v_1(t),\\
    x_n(t+1)=\min_{u\in\CU}\max_{w\in\CW}\{x_n(t)+\alpha_{uw} [x_{n-1}(t)-x_n(t)]+\beta_{uw}\}.
  \end{array}
\end{equation}
If we denote by $M^{uw}$ the matrices
$$M^{uw}=\begin{bmatrix}
        1 & 0 & \cdots & 0\\
        \alpha_{uw} & 1-\alpha_{uw} & & 0\\
        \vdots & \ddots & \ddots & \\
        0 & 0 & \alpha_{uw} & 1-\alpha_{uw}
      \end{bmatrix},$$
\vspace{5mm}
and by $c^{uw}(t)$ the vectors
$$c^{uw}(t)={}^t[v_1(t),\;\; \beta_{uw},\; \cdots,\; \beta_{uw},\;\; \beta_{uw}],$$
for $(u,w)\in\CU\times\CW$ and $t\in\BN$, then the dynamic system~(\ref{open1}) is written
\begin{equation}\label{open2}
  x_n(t+1)=\min_{u\in\CU}\max_{w\in\CW}\{[M^{uw} x(t)]_n+c^{uw}_n\},\quad 1\leq n\leq \nu\;.
\end{equation}

It is easy to check that the dynamic system~(\ref{open2}) is additive homogeneous of degree~1, and is
monotone under the condition $\forall (u,w)\in \CU\times\CW, \alpha_{uw} \in [0,1]$.
Therefore, under that condition,~(\ref{open2}) is non expansive.
However,~(\ref{open2}) is not connected for every $(u,w)\in \CU\times\CW$.

We will be interested here, in particular, in the stationary regime, where $v_1(t)$ reaches a fixed value $v_1$.
For this case, the eigenvalue problem associated to~(\ref{open1}) is
\begin{equation}\label{openeig}
  \begin{array}{l}
    \bar{v} + x_1=x_1+v_1,\\
    \bar{v} + x_n=\min_{u\in\CU}\max_{w\in\CW}\{x_n+\alpha_{uw} [x_{n-1}-x_n]+\beta_{uw}\}.
  \end{array}
\end{equation}
The system~(\ref{openeig}) is also written
\begin{equation}\label{open22}
  \bar{v}+x_n=\min_{u\in\CU}\max_{w\in\CW}\{[M^{uw} x]_n+c^{uw}_n\},\quad 1\leq n\leq \nu\;,
\end{equation}
where
$$c^{uw}={}^t[v_1,\;\; \beta_{uw},\; \cdots,\; \beta_{uw},\;\; \beta_{uw}].$$
Then we have the following result.
\begin{theorem}\label{openth1}
  For all $y\in\BR$ satisfying $\min_{u\in\CU}\max_{w\in\CW}(\alpha_{uw}y+\beta_{uw})=v_1$, the couple $(\bar{v},x)$
  is a solution for the system~(\ref{openeig}), where $\bar{v}=v_1$ and $x$ is given up to an
  additive constant by
  \begin{equation}\label{eqx}
    x={}^t[(n-1)y, \quad (n-2)y, \quad \cdots, \quad y, \quad 0].
  \end{equation}
\end{theorem}
\proof
  The proof is similar to that of Theorem~\ref{th-gen-2}.
  Let $y\in\BR$ satisfying $\min_{u\in\CU}\max_{w\in\CW}(\alpha_{uw}y+\beta_{uw})=v_1$.
  Let $(\bar{u},\bar{w})\in\CU\times\CW$ such that $\alpha_{\bar{u}\bar{w}}y+\beta_{\bar{u}\bar{w}}=v_1$.
  Let $x$ be given by~(\ref{eqx}). Then $\forall n\in\{1,2,\cdots,\nu\}$ we have
  \begin{align}
    [M^{\bar{u}\bar{w}} x]_n + c^{\bar{u}\bar{w}}_n & = (\alpha_{\bar{u}\bar{w}} y+\beta_{\bar{u}\bar{w}}) + x_n \nonumber \\
                  & = \min_{u\in\CU}\max_{w\in\CW}(\alpha_{uw} y+\beta_{uw}) + x_n \nonumber \\
                  & = \min_{u\in\CU}\max_{w\in\CW} [M^{uw} x]_n + c^{uw}_n \nonumber \\
                  & = v_1 + x_n. \nonumber
  \end{align}
\endproof

We can easily check that for $(u,w)\in\CU\times\CW$ such that $\alpha_{uw}=0$ and
$\beta_{uw}=v_1$, every inter-vehicular distance $y\in\BR$ satisfies the condition
$\min_{u\in\CU}\max_{w\in\CW}(\alpha_{uw}y+\beta_{uw})=v_1$. Thus, such couples $(u,w)$
do not count for that condition.
Theorem~(\ref{openth1}) can then be announced differently.
Let us denote by $\CW_u$ for $u\in\CU$ the family of index sets
$$\CW_u=\{w\in\CW, (\alpha_{uw},\beta_{uw})\neq (0,v_1)\},$$
and by $\tilde{y}$ the asymptotic average car-inter-vehicular distance:
\begin{equation}\label{diagop}
   \tilde{y}=\max_{u\in\CU}\min_{w\in\CW_u}\frac{v_1-\beta_{uw}}{\alpha_{uw}},
\end{equation}
where we use the convention $a/0=+\infty$ if $a>0$ and $a/0=-\infty$ if $a<0$.
Then Theorem~\ref{openth1} tells simply that if $\tilde{y}\in\BR$ (i.e. $-\infty<\tilde{y}<+\infty$), then
the dynamic system~(\ref{openeig}) admits a solution $(\bar{v},x)$ where
$\bar{v}=v_1$ is unique, and where $x$ is not necessarily unique and is given up to an additive constant by
$$x={}^t[(n-1)\tilde{y}, \quad (n-2)\tilde{y}, \quad \cdots, \quad \tilde{y}, \quad 0].$$

Non-uniform asymptotic car-distributions can also be obtained. Let us clarify the following three cases.
\begin{itemize}
  \item If $\exists u\in\CU, \forall w\in\CW_u, \alpha_{uw}=0$ and $\beta_{uw}<v_1$, then $\tilde{y}=+\infty$.
    In this case, the distance between the first car and the other cars increases over time and goes to $+\infty$.
    The asymptotic car-distribution on the road is not uniform.
  \item If $\forall u\in\CU, \exists w\in\CW_u, \alpha_{uw}=0$ and $\beta_{uw}>v_1$, then $\tilde{y}=-\infty$.
    In this case, the first car is passed by all other cars, and the distance between the first car and
    the other cars increases over time and goes to $+\infty$.
    The asymptotic car-distribution on the road is not uniform.
  \item If $\forall u\in\CU, \forall w\in\CW_u, \alpha_{uw}=0$ and if $\min_{u\in\CU}\max_{w\in\CW_u}\beta_{uw}=v_1$,
    then  for all $x\in\BR^{\nu}$, $(v_1,x)$ is a solution for the system~(\ref{openeig}).
    In this case, every distribution of the cars moving all with the constant velocity $v_1$ is stationary.
\end{itemize}

The formula~(\ref{diagop}) is the fundamental traffic
diagram expressing the average inter-vehicular distance as a function of the car
speed at the stationary regime. In the case where only
a \emph{minimum} operator is used in~(\ref{open1}), the formula~(\ref{diagop})
is reduced to the convex fundamental diagram
\begin{equation}\label{diagcx}
  y=\max_{u\in\CU}\frac{v_1-\beta_{u}}{\alpha_{u}}.
\end{equation}

\begin{example}

In order to understand the transient traffic, let us simulate the car-dynamics~(\ref{open1}).
We take as the time unit half a second (1/2 s), and as the distance unit 1 meter (m).
The parameters of the model are determined by approximating the behavior law~(\ref{Kerlaw}), with a free velocity
$v_0=14$ m/ 1/2s (which is about 100 km/h) and $\gamma=7.5$ as in~\cite{LWS99}; see Figure~\ref{app1}.

\begin{figure}[htbp]
  \begin{center}
    \includegraphics[width=10cm]{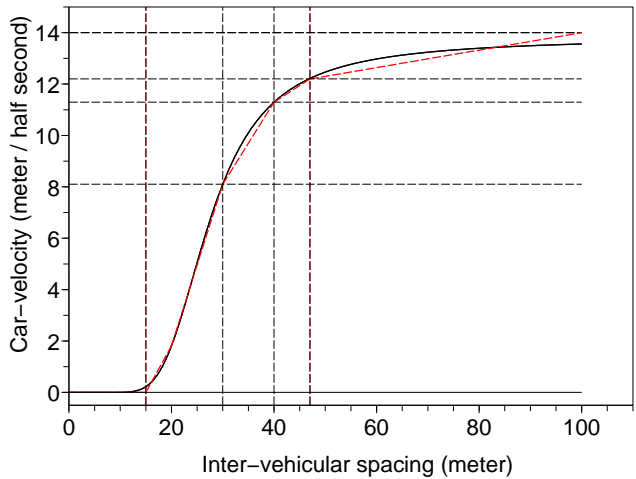}
    \caption{Approximation of the behavioral law~(\ref{Kerlaw}) with a piecewise linear curve.}
    \label{app1}
  \end{center}
\end{figure}

The behavior law is approximated by the following piecewise linear curve of six segments.
$$\tilde{V}(y)=\max\{\alpha_1 y+\beta_1, \min\{\alpha_2 y+\beta_2, \alpha_3 y+\beta_3, \alpha_4 y+\beta_4, \alpha_5 y+\beta_5,
       \alpha_6 y+\beta_6\}\},$$
where the parameters $\alpha_i$ and $\beta_i$ for $i=1,2, \cdots, 6$ are given by

 ~

\begin{tabular}{|c||c|c|c|c|c|c|}
  \hline
  Segments & 1 & 2 & 3 & 4 & 5 & 6 \\
  \hline
  $\alpha_i$ & 0 & 0.54 & 0.32 & 0.13 & 0.34 & 0 \\
  \hline
  $\beta_i$ & 0 & -8.1 & -1.47 & 6.11 & 10.6 & 14\\
  \hline
\end{tabular}

~

We simulate the piecewise linear car-following model associated to the approximation above.
\begin{equation}
  \begin{array}{l}
     x_1(t)=x_1(t-1)+v_1(t), \\
     x_n(t)=x_n(t-1)+\tilde{V}(x_{n-1}(t-1)-x_n(t-1)).
  \end{array}
\end{equation}
The velocity of the first car $v_1(t), t\geq 0$ is varied in the time interval $[0,1000]$, then
fixed to the free velocity $v_0=14$ m/ 1/2s in the time interval $[1000,3000]$, and finally fixed on
a velocity that exceeds $v_0$ in the remaining time $[3000,7200]$. The average inter-vehicular distance
is then computed at every time $t$, and the results are shown in Figure~\ref{vel-inter2}.

\begin{figure}[htbp]
  \begin{center}
    \includegraphics[width=8.5cm]{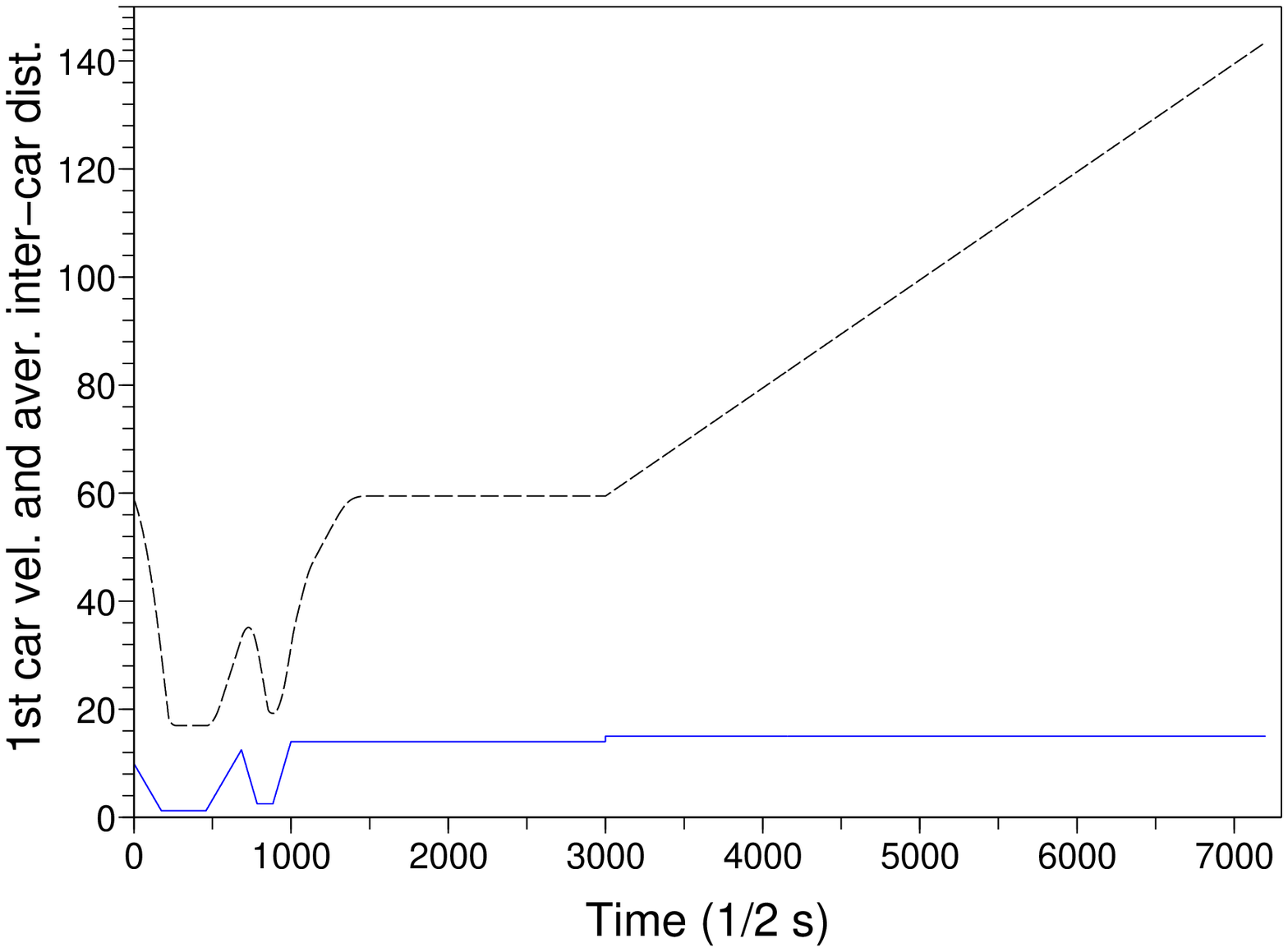} \hspace{-1cm}
    \includegraphics[width=8.5cm]{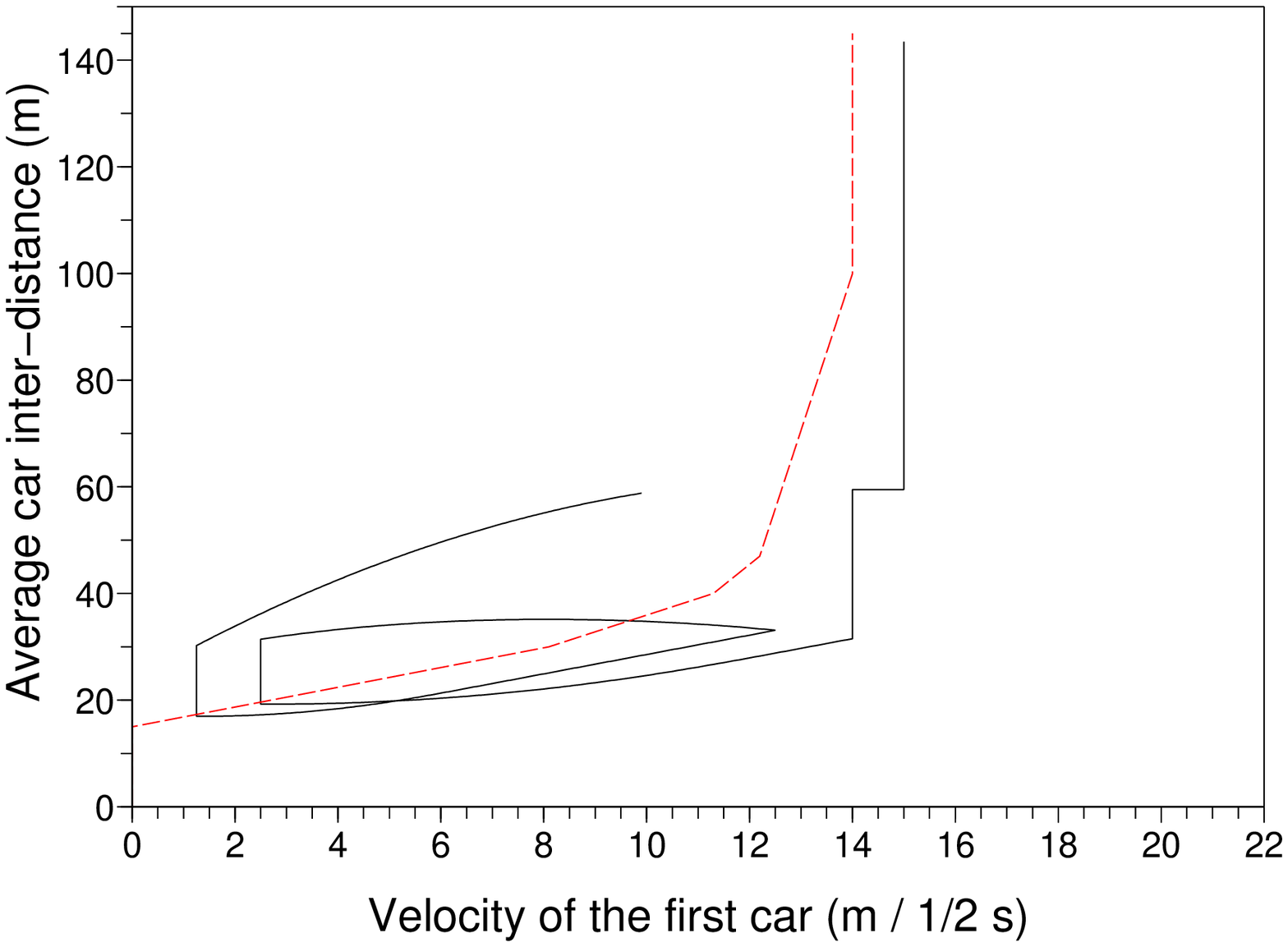}
    \caption{Simulation results. On the left-side: the first car velocity (solid line), and the average inter-vehicular distance
       (dash line) functions of time. On the right side: the approximation of the behavior law~(\ref{Kerlaw}) (dash line), and
       the average inter-vehicular distance obtained by simulation in function of the velocity of the first car (solid line).}
    \label{vel-inter2}
  \end{center}
\end{figure}

The simple simulation we made here permits to have an idea of the traffic in the transient regime.
Figure~\ref{vel-inter2} shows how the average of the inter-vehicular distance is changed due to a changing in the
velocity of the first car. In the right side of Figure~\ref{vel-inter2}, we compare the fundamental law assumed
in the model with the diagram giving the \emph{average} inter-vehicular distance (with respect to the number of cars)
function of the velocity of the first car (a kind of macroscopic fundamental diagram).
We observe that loops are obtained on that diagram in the transient traffic.
The loops are interpreted by the fact that once the velocity of the first car is temporarily stationary, the
velocities of the following cars, and thus also the average velocity of the cars, make some time to attain
the first car velocity. It can also be interpreted by saying that even though the cars have, individually,
the same response to a changing in inter-vehicular distance; their collective response depends on whether the
inter-vehicular distance is increasing or decreasing.
The apparition of such loops is due to the reaction time of drivers.
It can also be related to the number of anticipation cars in case of multi-anticipative modeling.
However, one may measure on a given section, different car-flows for the same car density (or occupancy rate)
depending on the traffic acceleration or deceleration, and interpret it as the hysteresis
phenomenon~\cite{Edi61,TM74,Zha99,GS11}.

\end{example}

\section{Conclusion}

We proposed in this article a car-following model which extends
the linear car-following model as well as the min-plus model.
The stability and the stationary regimes of the model are characterized
thanks to a variational formulation of the car-dynamics.
This formulation, although already made with continuous-time models,
it permits to clarify the stimulus-response process in microscopic
discrete-time traffic models, and to interpret it in term of
stochastic optimal control. Among the important questions to be treated in the
future, the impacts of heterogeneity and anticipation in driving,
on the transient and stationary traffic regimes, based on the model proposed
in this article.

\appendix

\section{Connectedness of system~(\ref{gener1})}
\label{appA}

Let $x\in\BR^{\nu}$. We denote by $h: \BR^{\nu} \to \BR^{\nu}$ the operator defined by $h(x(t))=x(t+1)$,
where $x_n(t+1)$, for $1\leq n\leq \nu$ are given by the definition of the system~(\ref{gener1}).
That is
\begin{equation}\nonumber
  x_n(t+1) = \min_{u\in\CU}\{[M^u x(t)]_n+c^u_n\},\quad 1\leq n\leq \nu\;.
\end{equation}

\begin{itemize}
  \item If $\exists u\in\CU$, such that $\alpha_u\in(0,1]$, then for all $n\in\{1, 2, \cdots, \nu\}$,
    there exists an arc, on the graph associated to $h$, going from $n-1$ to $n$
    ($n$ being cyclic in $\{1, 2, \cdots, \nu\}$).
    Indeed,
    $$x_n(t+1)=(1-\alpha_u)x_n(t)+\alpha_u x_{n-1}(t)+\beta_u.$$
    Then since $\alpha_u>0$, we get:
    $$\lim_{\nu\to\infty}h_n(\nu e_{n-1})=\lim_{\nu\to\infty}[\alpha_u \nu+\beta_u]=\infty.$$
    where $e_{n-1}$ denotes the $(n-1)^{th}$ vector of the canonic basis of $\BR^{\nu}$.
    Therefore the graph associated to $h$ is strongly connected.
  \item If $\forall u\in\CU, \alpha_u=0$, then we can easily check
    that all arcs of the graph associated to $h$ are loops. Hence that graph is not strongly
    connected.
\end{itemize}

\section{Duality in traffic modeling}
\label{appB}

We show in this appendix the duality in traffic modeling of using the three functions
\begin{itemize}
  \item $n(t,x)$: cumulated number of cars passed through position $x$ from time $0$ up to time $t$.
  \item $x(n,t)$: position of car $n$ at time $t$ (or cumulated traveled distance of car $n$ from time
                       $0$ up to time $t$).
  \item $t(n,x)$: the time that car $n$ passes by position $x$.
\end{itemize}
We base on the Lagrangian traffic descriptions given in the introduction,
and give the equivalent traffic descriptions by using the functions $n(t,x)$ and $t(n,x)$,
in both cases of discrete-time and continuous-time modeling.



\begin{enumerate}


\item In Eulerian traffic descriptions, the function $n(t,x)$ is used.
The partial derivative $\partial_t n(t,x)$ expresses the car-flow $q(t,x)$ at time $t$ and position $x$,
while $- \partial_x n(t,x)$ expresses the car-density $\rho(t,x)$ at time $t$ and position~$x$.
The equality $\partial_{tx} n(t,x)=\partial_{xt} n(t,x)$ gives the car conservative law:
\begin{equation}\label{conservlaw}
  \partial_t k (t,x) + \partial_x q(t,x) = 0.
\end{equation}
The first order traffic model LWR~\cite{LW55,Ric56} supposes the existence of a fundamental diagram of traffic
giving the car-flow $q$ as a function of the car-density $\rho$ at the stationary regime, through a
function $Q_e$, and that the diagram also holds for the transient traffic~:
\begin{equation}\label{behavlaw}
  q(t,x)=Q_e(\rho(t,x)).
\end{equation}
Then~(\ref{conservlaw}) and~(\ref{behavlaw}) give the well known LWR model
\begin{equation}\label{LWR}
  \partial_t \rho (t,x) + \partial_x \rho(t,x) Q_e'(\rho) = 0.
\end{equation}

Also discrete-time and -space Eulerian traffic models exist. Those models are
in general derived from Petri nets as in~\cite{BCOQ92,Far08,FGQ11}. For example the
traffic on a 1-lane road without passing can be described by
\begin{equation}\label{dynntx}
  n(t,x)=\min\{a_x + n(t-\tau_x, x-\delta x), \quad \bar{a}_{x+\delta x} + n(t-\bar{\tau}_{x+\delta x}, x+\delta x)\},
\end{equation}
where
\begin{itemize}
  \item[$a_x$] denotes the number of cars being in $(x-\delta x,x)$ at time zero.
  \item[$\tau_x$] denotes the free travel-time of a car from $x-\delta x$ to $x$.
  \item[$\bar{a}_x$] denotes the space non occupied by cars in $(x-\delta x,x)$ at time zero.
         If we denote by $c_x$ the maximum number of cars that can be in $(x-\delta x,x)$, then
         we have simply $\bar{a}_x=c_x-a_x$.
  \item[$\bar{\tau}_x$] denotes the reaction time of drivers in $(x-\delta x,x)$. That is the time
         interval between the time when $(x-\delta x,x)$ is free of cars and the time when a car
         being in $(x-2\delta x,x-\delta x)$ starts moving to $(x-\delta x,x)$.
         If we denote by $T_x$ the total traveling time (reaction time + moving time) of a car
         from $(x-2\delta x,x-\delta x)$ to $(x-\delta x,x)$, then we have simply
         $\bar{\tau}_x = T_x - \tau_x$.
\end{itemize}


\item By using the variable $t(n,x)$, the first order differentiation of $t(n,x)$ with respect to $n$,
denoted $z$ is $z(n,x)=-(t(n-1,x)-t(n,x))$, while the derivative of $t(n,x)$ in $x$, denoted $r$ is
$r(n,x)=\partial_x t(n,x)$. We notice here that $z$ and $r$ are interpreted respectively as the inverse
flow and the inverse velocity of vehicles. A conservation law (of time) is then written
\begin{equation}\label{conservlaw3}
  \partial_x z(n,x)+r(n,x)-r(n-1,x) = 0.
\end{equation}
The law (\ref{conservlaw3}) combined with the fundamental diagram $r=R_e(z)$ gives the model
\begin{equation}
  \partial_x r(n,x) = R'(z(n,x)) \Delta r(n,x),
\end{equation}
where $\Delta r(n,x)=r(n-1,x)-r(n,x)$.
Note that, having a fundamental diagram $v=V_e(q)$ giving the stationary velocity as a function of
the stationary flow, the diagram $R_e$ is nothing but $R_e(z)=1/V_e(1/z)$.

Discrete-time-and-space modeling with the function $t(n,x)$ also exist.
The models are also inspired from Petri net, and dual dynamics to~(\ref{dynntx}) are obtained.
For example, using the same notations as in~(\ref{dynntx}), the traffic on a 1-lane road without
passing can be described by
\begin{equation}\label{dyntnx}
  t(n,x)=\max\{\tau_x + t(n-a_x, x-\delta x), \; \bar{\tau}_{x+\delta x} + t(n+\bar{a}_{x+\delta x}, x+\delta x)\}.
\end{equation}
Note here that a \emph{max} operator is used rather than a \emph{min} one.
For more details on the duality of~(\ref{dynntx}) and~(\ref{dyntnx}) and the meanings in term
of Petri nets, event graphs and min-plus or max-plus algebras, see~\cite{BCOQ92}.

\end{enumerate}


\bibliographystyle{plain}
\bibliography{biblio}

\end{document}